\documentclass[a4paper, 11pt]{article}
\usepackage{amsmath}
\usepackage{amsfonts}
\usepackage{textcomp}
\usepackage{mathcomp}
\usepackage{lmodern}
\usepackage{anyfontsize}
\usepackage{Zoe_C_IASS_Style}
\usepackage{color,ifpdf,latexsym,url}

\usepackage  [breaklinks,bookmarks,bookmarksnumbered,bookmarksopen,bookmarksopenlevel=2]{hyperref}
{\makeatletter \hypersetup{pdftitle={\@title}}}
\def\authornames{Zoe COOPERBAND, Allan McROBIE, Cameron MILLAR, Bernd SCHULZE} %HERE, Please replace "the authors" with actual author names like the following:

\usepackage[sorting=none,style=ieee,citestyle=numeric]{biblatex}
\newcommand{\tblcaption}[1]{\def\@captype{table}\caption{#1}}

\begin{document}

\vspace*{0mm}
\begin{center}%

{\LARGE \bfseries Homology of Moment Frames\\}%
\vspace{1.0em}
{\normalsize Zoe COOPERBAND$^\textrm{*}$, Allan McROBIE$^a$, Cameron MILLAR$^b$, Bernd SCHULZE$^c$, \\}
\vspace{1.0em}
{\small $^*$ Department of Electrical and Systems Engineering; GRASP Lab\\University of Pennsylvania\\zcooperband\{at\}gmail\{dot\}com\\}
\vspace{1em}
{\small $^\textrm{a}$ Cambridge University Engineering Dept, Cambridge, CB2 1PZ, UK \\} 
{\small $^\textrm{b}$ Skidmore Owings \& Merrill, The Broadgate Tower, 20 Primrose St, London, EC2A 2EW, UK \\}
{\small $^\textrm{c}$  Dept of Mathematics \& Statistics, Lancaster Univ., UK }
\vspace{1.0em}

\end{center}%
%\vspace{1.0em}
\begin{abstract}
Using homological techniques 
%from applied algebraic topology are used to develop a novel conceptual framework that encompasses the fundamental structural properties of both axially-loaded pin-jointed trusses and rigid moment frames. 
%We focus on the use of cellular cosheaves, which enable high level analysis and decomposition of the geometric data that describes the kinematic and static properties of frames. We 
we show that a pin-anchored frame that involves only moments and shears provides a conceptual bridge between the statics of moment frames and the kinematics of pin-jointed trusses. One immediate result is a long exact sequence whose alternating sum of dimensions gives a novel counting rule for self-stresses and mechanisms. This combines the Maxwell-Calladine count for pin-jointed trusses with the circuit rank (first Betti number) associated with self-stresses in moment frames. These relations apply to frames in 2, 3 or any dimensions. This work heralds a shift towards a deeper study of the relationships and dualities that exist between structural equilibria and kinematics.%\textcolor{red}{I think we should shorten this! (e.g. the first two sentences can basically go. Third setnence can start with: Using homological techniques we show.}
\end{abstract}
  
\vskip 1.0em

\section{Introduction}
%Homological techniques from applied algebraic topology provide a conceptual framework in which to describe the fundamental structural properties of both axially-loaded pin-jointed trusses and rigid moment frames.
Homological algebra provides powerful tools for describing the fundamental structural properties of both axially-loaded pin-jointed trusses  and rigid moment frames~\cite{TowardsCooperband2023,Part1Cooperband2023, Zoe2024equivariant}. This paper focuses on the use of {\it cellular cosheaves}~\cite{Curry2013}. These are used to store the geometric data embodied by the kinematic and static properties of trusses and frames:  the forces, moments, displacements and rotations and their inter-relationships that are familiar concepts in structural engineering. Once the data is so organised, the full power of the description arises when the techniques of homology theory are applied. Specifically, counting rules that relate node and bar counts to the numbers of possible mechanisms and states of self-stress are developed, creating a novel link between two familiar counts:  
\begin{itemize}
    \item the Maxwell-Calladine count for a {\it pin-jointed  truss}~\cite{Calladine1978}; and 
    \item the (circuit rank) count of the number of cuts needed to make a statically-determinate tree from a {\it rigidly-connected moment frame}.
\end{itemize}

Central to this connection is an intermediate frame system whose joints and members can transmit shears and moments, but whose bars cannot carry axial tension. We call this an ``anchored frame.''

This work represents a first step towards a rigidity theory for moment-resisting frames. Following pioneering work by  Henneberg~\cite{henneberg1903}, Pollaczek-Geiringer \cite{PollaczekGeiringer} and Laman~\cite{laman1970} among others, much 
is now known about the rigidity of structures. Both Graver's classic text {\it Counting on Frameworks}~\cite{graver2001} and the recent and more comprehensive interdisciplinary book {\it Frameworks, Tensegrities, and Symmetry} \cite{connelly_guest_2022} by  Connelly and  Guest not only provide strong introductions but also highlight an immediate problem: in their study of geometric %constraint
systems with distance constraints, rigidity theorists concern themselves almost exclusively with what they call {\it bar-joint frameworks}. In structural engineering, these are referred to as {\it pin-jointed trusses}, with bars transmitting only axial loads. In this paper, we widen the focus to the rigidity of a wider class of structures whose inter-connected bars can carry not only axial forces but also shear forces, moments and torsions. In engineering, it is this ability to transmit moments that warrants the use of the word {\it frame} or {\it framework}.  

This paper is thus a step towards a possible substantial broadening of scope for rigidity theory, opening up potential avenues for the detailed study of the rigidity of {\it frames}. Although arguments exist about why axial-only structures may offer material efficiencies, beams that bend are central to structural engineering and much of the infrastructure around us.  Moment frames are thus worthy of deeper study. 

Beyond this, a theory that can sensibly represent moments may also help resolve certain singularity issues that arise within the pin-jointed, axial-only ansatz. There, it is well known that some topologies can carry states of self-stress only if certain geometrical conditions are satisfied. A classic example is a 2D truss in the Desargues configuration (see Fig IV of Maxwell 1864 \cite{maxwell1864a} and Fig.~\ref{fig:moment_desargues} later) where the truss can carry a  state of self-stress (and possess a mechanism) only if the lines of three particular bars meet at a point. This is a highly non-generic condition. In practice, just how close do those three bars have to be to meeting at a single point? If moments are considered this all-or-nothing pathology is avoided. If the three bars do not quite meet at a point, states of self-stress still exist nearby: they will simply contain some small moments. The anchored frame resolves some of these concerns, describing where truss self-stresses and mechanisms algebraically ``go'' when shifting the geometry.

\section{Cosheaves and Truss Statics}

We develop frame and truss statics in terms of cellular cosheaves and their homology. Although this is a highly abstract definition, one quickly finds that cosheaves simply redescribe systems and operations that are intuitive and used by engineers everyday.

A {\it framework} is an abstract graph $G = (V,E)$ along with a geometric realization $p: V \to \R^n$  %or $p: V \to \R^3$
that assigns each vertex $v\in V$ a %geometric coordinate
point $p(v) = p_v$ in Euclidean space. Edges are realized as straight line segments with %the pair
$(G, p)$ a geometric graph. Whenever vertices $u,v$ are incident to an edge $e$ we use the notation $u,v\lhd e$. In this paper we will primarily consider frameworks in the ambient space $\R^2$.

%For notation purposes let $p_e = \frac{p_v + p_u}{2}$ denote the center of an edge. 

% These edges are assigned an arbitrary orientation with a {\it signed incidence function} detecting this orientation; using the notation $v,u\lhd e$ for two vertices incident to the edge $e$, we say $[v:e] = +1$ if $e$ is oriented towards $v$ and $[u:e] = -1$ if $e$ is orineted away from $u$.

\begin{definition}[Cellular cosheaf,~\cite{Curry2013}]
    Over a framework $(G, p)$ a {\it cellular cosheaf} $\calk$ is comprised of
    \begin{itemize}
        \item vector spaces $\calk_e$ and $\calk_v$ assigned to each edge $e$ and vertex $v$, respectively, called {\it stalks}, and
        \item linear maps $\calk_{e\rhd v}: \calk_e \to \calk_v$ between every incidence $v\lhd e$.
    \end{itemize}
\end{definition}
A cellular cosheaf as of now is nothing more than a way to bookkeep what geometric data is assigned to which cell, and how this data interacts.

\begin{example}[2D truss statics]
    What is the stress data assigned to a pin-jointed truss? Each edge $e$ of a framework $(G, p)$ is assigned an axial tensile or compressive force %, a real-valued
    scalar. These internal forces are propagated to vertices, where these forces are summed with external loads. The axial {\it force cosheaf} encoding these forces has $\calf_e = \R$, $\calf_v = \R^2$, and linear maps $\calf_{e\rhd v}$ being embeddings of $\R$ into the larger space $\R^2$. In other words, at an edge $u,v\lhd e$ both maps $\calf_{e\rhd v}$ and $\calf_{e\rhd u}$ are the same $2\times 1$ %block
    matrix $\ell_e = (p_{uy} - p_{vy}, p_{ux} - p_{vx})^\top$, with range denoted by the subspace ${\bf e} \subset \R^2$. %parallel with the embedding of $e$.

    The truss {\it equilibrium matrix} ${\bf A}$ is a size $2|V| \times |E|$ matrix that transmits all internal loads to external load spaces, where $|\cdot|$ denotes the dimension/size of the space/set. (In  rigidity theory this matrix is the transpose of the so-called \emph{rigidity matrix} \cite{SWbook}.) If $w$ is a vector of axial loads over edges, we say that $w$ is a {\it self stress} if ${\bf A}(w) = 0$. This condition requires all internal forces to equilibrate at each vertex without the need of external forces. In cosheaf notation, we combine all linear maps into the {\it boundary matrix} $\bdd_\calf = {\bf A}$ with non-zero blocks $\pm\calf_{e\rhd v}$ wherever there is an incidence $v\lhd e$. The axial load assignment $w$ is called a {\it chain} with local components $w_e\in \calf_e = \R$. If $\bdd_\calf (w) = 0$ then $w$ is called a {\it cycle}.
\end{example}

The terminology of chains and cycles comes from algebraic topology. Chains are arbitrary vector assignments to all stalks of a given dimension, that is, elements of the combined stalk spaces
\begin{equation}
    C_0 \calk = \bigoplus_{v\in V} \calk_v; \qquad C_1 \calk = \bigoplus_{e\in E} \calk_e,
\end{equation}
(the same as the product $\calk_{v_1} \times \calk_{v_2} \times \cdots$ of finite dimensional vector space stalks). Furthermore, the cosheaf boundary map is a linear map between spaces of chains
\begin{equation}
\label{eq:bdd_cosheaf}
    \bdd_\calk: C_1 \calk \to C_0 \calk; \qquad (\bdd_\calk(w))_v = \sum_{e: v\lhd e} \pm \calk_{e\rhd v}(w_e),
\end{equation}
where $(\bdd_\calk(w))_v$ is the component of $\bdd_\calk(w)\in C_0 \calk$ in the stalk $\calk_v$. The choice of sign in Equation~(\ref{eq:bdd_cosheaf}) is informed by an arbitrary orientation for the edges -- the sign is positive if the edge ``points towards'' $v$ and is negative otherwise. The {\it homology} of a cosheaf $\calk$ %over a framework $(G, p)$
then consists of the two vector spaces
\begin{equation}
    H_1 \calk = \ker \bdd_\calk \subset C_1 \calk; \qquad H_0 \calk = C_0 \calk / \im \bdd_\calk \iso (\im \bdd_\calk)^\perp \subset C_0 \calk
\end{equation}
where the latter space is a quotient vector space of $C_0 \calk$ by the image of the boundary map $\bdd_\calk$.

Homology often captures the most important aspects of the system; for the force cosheaf, $H_1 \calf$ is the space of self stresses while $H_0 \calf$ is (isomorphic to) the space of infinitesimal degrees of freedom, the combined space of rigid body DOF $\scrr \calf$ and mechanisms $\scrm \calf$.

\begin{example}[Classical homology]
    Chain complexes and homology were born from classical constructions in algebraic topology from a century ago~\cite{AlgebraicHatcher2002}. Invented to quantify topological features, classical homology detects the number and arrangement of voids and holes of an abstract space. One attaches a scalar value $\R$ to every cell, with boundary map $\bdd: C_i G \to C_{i-1} G$ called the {\it signed incidence matrix}. Over a graph $G$, $H_1 G$ is the linear span of cycles while $H_0 G$ encodes connected components. %This construction can be modeled by a {\it constant cosheaf} $\overline{\R}$ if desired~\cite{Part1Cooperband2023}.

    The classical {\it Euler characteristic} $\mathfrak{X}$, the alternating sum of vector space dimensions, is a powerful topological invariant of spaces~\cite{AlgebraicHatcher2002}. For a graph $G$, using the rank-nullity theorem we find
    \begin{equation}
    \label{eq:euler}
        \mathfrak{X}(C G) = |C_0 G| - |C_1 G| = |V| - |E| = |H_0 G| - |H_1 G| = \mathfrak{X}(H G).
    \end{equation}
    If $G$ is connected then $H_0 G$ is $1$-dimensional and following Equation~\eqref{eq:euler} the number of graph cycles is $|H_1 G| = |E| - |V| + 1$ (in graph theory, this is known as the {\it circuit rank}).
\end{example}

An example of the Euler count applied to cosheaves is the {\it Maxwell-Calladine counting rule}~\cite{Calladine1978}
\begin{equation}
\label{eq:maxwell_rule}
    \mathfrak{X}(C \calf) = 2|V| - |E| = |\scrr \calf| + |\scrm \calf| - |H_1 \calf| = \mathfrak{X}(H \calf)
\end{equation}
over the force cosheaf $\calf$, enveloping counts of self-stresses, rigid DOF $\scrr \calf$, and mechanisms $\scrm \calf$~\cite{TowardsCooperband2023}. This is the alternating count of Equation~\eqref{eq:euler} replacing standard homology with the homology of $\calf$.

\section{Frame Statics}

\begin{figure}[ht]
    \centering
    \begin{subfigure}[t]{0.49\textwidth}\centering
        \includegraphics[scale = 1.2]{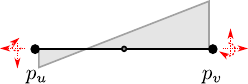}
        \caption{The moment cosheaf $\calm$ over an edge. The rigid beam transmits axial, shear, and bending forces. A moment function is graphed over the beam.}
    \end{subfigure}
    \hfill
    \begin{subfigure}[t]{0.49\textwidth}\centering
        \includegraphics[scale = 1.2]{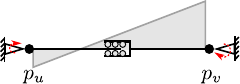}
        \caption{The anchored cosheaf $\caln$ over an edge. The mid-beam joint transmits bending, shear, and torque but no axial force. The pinned anchors remove residual shear.}
    \end{subfigure}
    \vspace{1pt}
    \caption{Sketches of the moment $\calm$ and anchored $\caln$ cosheaves over an edge.}
    \label{fig:moment_anchored}
\end{figure}

\begin{example}[Plane frame statics]
    Moment frame elements have a simple formulation and boundary map. Suppose that $(G, p)$ is a single edge with two endpoints $u,v\lhd e$ embedded horizontally along the $x$-axis in $\R^2$. Letting $M$ be a function representing the moment along the beam; it is well known that its derivative $\frac{dM}{dx}$ is the shear along the beam. Assume that the shear force $F_y$ is constant so $M$ has constant slope. Then the linear function $M$ over the edge $e$ is $M(t) = M_e + F_y \cdot t$, where $M_e$ is the moment at the edge center $\frac{p_u + p_v}{2}$ and $t$ parameterizes a location along the length of the beam; see Fig.~\ref{fig:moment_anchored} (a) where a moment function $M(t)$ is graphed over an edge. Then the induced moments at $p_u$ and $p_v$ are the magnitudes of $F_y \frac{p_u - p_v}{2}$ and $F_y \frac{p_v - p_u}{2}$. Letting $w_e = (M, F_x, F_y)$ be a force-couple at the edge center, the induced spatial force at $p_u$ and $p_v$ in coordinates follows from applying matrices:
    \begin{small}\begin{equation}
    \label{eq:moment_boundary}
        \bdd_\calm(w_e)_u =
            \begin{bmatrix}
                 1 &  0 & \|\frac{p_u - p_v}{2}\| \\
                 0 & 1 & 0 \\
                 0 & 0 & 1
            \end{bmatrix}
            \begin{bmatrix}
                 M \\
                 F_x \\
                 F_y
            \end{bmatrix}
        ;
        \qquad
        \bdd_\calm(w_e)_v = -
            \begin{bmatrix}
                 1 &  0 & -\|\frac{p_u - p_v}{2}\| \\
                 0 & 1 & 0 \\
                 0 & 0 & 1
            \end{bmatrix}
        \begin{bmatrix}
                 M \\
                 F_x \\
                 F_y
            \end{bmatrix}
    \end{equation}\end{small}
    where the edge $e$ is oriented towards $u$ in this notation.
\end{example}

To cleanly represent a force couple transform, we introduce the {\it exterior product space} $\bigwedge^2 \R^n$, the vector space consisting of formal bilinear pairs $x\wedge y$ where $x\wedge y = -y \wedge x$ over vectors $x,y\in \R^n$. It follows immediately from this definition that $x\wedge x = 0$. With $x,y$ a basis for $\R^2$ the pseudo-scalar $M = x\wedge y$ is a basis for $\bigwedge^2 \R^2$. In dimension three, $\bigwedge^2 \R^3$ consists of moments $M_x, M_y, M_z$ in the three directions, with wedge product equivalent to the cross product (after applying the Hodge-star operator).

The exterior product models moments. The moment generated by a force vector $F\in \R^n$ applied at a lever arm $\ell\in \R^n$ is equal to the product $F\wedge \ell$.

\begin{definition}[Moment cosheaf]
\label{def:cosheaf_moment}
    Over a framework $(G, p)$ in $\R^2$, the {\it moment cosheaf} $\calm$ has stalks $\calm_e = \calm_v = \bigwedge^2\R^2 \oplus \R^2$ comprised of force-couples at each cell. Each cosheaf stalk map
    \begin{equation}
    \label{eq:moment_extension}
        \calm_{e\rhd v}(M, F) = (M + F\wedge \frac{p_v - p_u}{2}, F)\in \calm_v
    \end{equation}
    sends a force couple at the edge center $\frac{p_u + p_v}{2}$ to a force couple at the coordinate $p_v$.
\end{definition}

One can check that the linear map~\eqref{eq:moment_extension} aligns with those in Equation~\eqref{eq:moment_boundary}, combining to form the size $3|V| \times 3|E|$ frame equilibrium matrix ${\bf B} = \bdd_\calm: C_1 \calm \to C_0 \calm$ after choosing bases. The moment cosheaf $\calm$ has a particularly simple homology type. In a connected rigid frame clearly the combined system has three rigid body DOF in $\R^2$ and six rigid body DOF in $\R^3$, spanning $H_0 \calm = \scrr \calm$.

The homologies $H_1 \calm$ and $H_0 \calm$ are in fact isomorphic to $3$ copies of the classical homologies $H_1 G$ and $H_0 G$ in $\R^2$ and $6$ copies in $\R^3$. Applying the Euler equation~\eqref{eq:euler} to the cosheaf $\calm$ in $\R^2$ we find 
\begin{equation}
\label{eq:euler_moment}
    |H_1 \calm| = |H_0 \calm| + 3( |E| - |V|) = 3 ( |E| - |V| + 1) = 3|H_1 G|.
\end{equation}
where $H_1 \calm = \ker {\bf B}$ is the space of frame self stresses. %\textcolor{red}{Is there a reference for this?}
This circuit rank count is familiar in elementary structural engineering. Trivially, a fully-welded frame requires $|E| - |V| + 1$ cuts to make it a statically-determinate tree.  In 2D there are three independent stress resultants (axial, shear and bending) at each such cut, hence Equation~(\ref{eq:euler_moment}). In 3D there are six stress resultants at each cut.

\section{The Anchored Cosheaf}

While the axial force cosheaf $\calf$ and moment cosheaf $\calm$ have been discussed  separately, the true power of cosheaf theory comes from the relationship between the two.

\begin{definition}[Cosheaf map]
\label{def:cosheaf_map}
    A {\it map between cosheaves} $\phi: \calk \to \call$ is a collection of stalk-wise linear maps $\phi_e: \calk_e \to \call_e$ and $\phi_v: \calk_v \to \call_v$ such that the map compositions $\call_{e\rhd v} \phi_e = \phi_v \calk_{e\rhd v}$ are equal at every incidence $v\lhd e$. This equation ensures that the systems ``align/agree'' across cell incidences.
\end{definition}

We can formally state how the moment cosheaf $\calm$ subsumes the force cosheaf $\calf$, meaning truss stresses ``are'' frame stresses. There is an injective map $\phi: \calf \to \calm$ which is comprised of embeddings $\phi_e:\calf_e \to 0 \oplus \R^2 \subset \calm_e$ at edges. At vertices $\phi_v: \calf_v \to 0 \oplus \R^2 \subset \calm_v$ simply adds zero moments to nodal forces in $\calf_v$. To see this satisfies Definition~\ref{def:cosheaf_map}, if $\phi_e(w) = F_e\in {\bf e} \subset \R^2$ is an axially aligned force vector then $F_e \wedge \frac{p_u - p_v}{2} = 0$ in the stalk map $\calm_{e\rhd v}$, as $(p_u - p_v)\in {\bf e}$. Thus there is no moment component and $\phi_v \calf_{e\rhd v} = \calm_{e\rhd v}\phi_e$.

The injective cosheaf map $\phi$ induces a quotient cosheaf $\calm / \phi\calf$ which we will denote as $\caln$ and call a {\it anchored cosheaf}. This quotient cosheaf has stalks $\caln_e = \calm_e / \phi\calf_e = \bigwedge^2 \R^2 \oplus \R^2/ {\bf e}$ of dimension $2$ and stalks $\caln_v \iso \bigwedge^2 \R^2$ of dimension $1$ over edges and vertices. At the physical level, the cosheaf $\caln$ ignores the axial forces over edges and ignores all forces at vertices, so that it only considers the moment and shear components. A cosheaf stalk map $\caln_{e \rhd v}$ takes an axial-free force couple in $\caln_e$ and computes its induced moments at the edge endpoints, quotienting out the forces at $v$.

\begin{example}[Statics of the anchored cosheaf]
\label{ex:anchored_frame}
    An {\it anchored frame} should be thought of as a moment frame where each stiff juncture $v$ is held in place by an external anchor restricting translation but not rotation. These ``anchors'' are physically pin-joints to an external system\footnote{It is critical to note that the members are still rigidly attached to one another: the junction is not a pin itself but instead has reaction forces applied to it through a pinned connection.}. However, there then are are trivial degrees of axial self-stress over each edge; to eliminate these we insert a prismatic ``sliding joint'' mid-member that permits free axial extensions but restricts (and transmits) shear and moments. A sketch of an anchored frame element is pictured in Fig.~\ref{fig:moment_anchored} (b). The extending joint has been previously utilized to model shell structures~\cite{calladine1983}.

    The anchored cosheaf $\caln$ over a framework $(G, p)$ models this anchored frame system. The edge stalk $\caln_e \iso \bigwedge^2 \R^2 \oplus \R^2/ {\bf e}$ has the axial force quotiented out -- equivalent to the mid-member extension joint. The vertex stalks $\caln_v \iso \bigwedge^2 \R^2$ only need to detect moments because the anchor, as a pin-joint, absorbs any residual force. The homology $H_1 \caln$, the kernel of the size $|V| \times 2|E|$ equilibrium matrix $\bdd_\caln$, describes the valid states of self-stress of the anchored system. Moments and shears combine just as in a moment frame (in $\calm$), but the external pinned reactions eliminate the residual net shear force.
\end{example}

% \begin{remark}[``Non-degenerate'' Anchor Frame]
% \label{rem:non-degenerate_anchor}
%     For any ``non-degenerate'' geometric graphs $(X^1, o)$, the homology $H_0 \caln$ will always be trivial. By non-degenerate we mean that $X$ each vertex $v$ has incident edges $e_i$ such that the collection of image subspaces $\im \calk_{e_i \rhd v} = \R^n \wedge {\bf e_i}$ combine to span $\bigwedge^2 \R^n$. In two dimensions because $\bigwedge^2 \R^2 = \R^2 \wedge {\bf e_i}$ at any edge $e_i$, this condition simply means that $X$ is not a point. In $\R^3$ this means that every vertex $v$ has at least two incident and non-collinear edges $e_i, e_j$: the subspace $\R^3 \wedge {\bf e_i}$ associated to every edge $e_i$ is two dimensional, lacking the bivector encoding torsion along the member. This non-degeneracy assumption eliminates potential cases of the ``spinning node'' problem that often arises in computational statics.
% \end{remark}

In summary, the cosheaves $\calf, \calm,\caln$ form a short exact sequence of cosheaves
\begin{equation}
\label{eq:force_moment_ses}
    0 \to \calf \xrightarrow{\phi} \calm \xrightarrow{\pi} \caln \to 0.
\end{equation}
where $\pi: \calm \to \caln$ is a cosheaf map projection. Exactness means that $\ker \pi = \im \phi$; consequently stalks can be decomposed as $\calm_e\iso \calf_e \oplus \caln_e$ and $\calm_v \iso \calf_v \oplus \caln_v$. The map $\phi$ converts a pin-jointed truss to a moment frame by physically ``gluing'' the pins shut, allowing them to transfer shear and moments. The map $\pi$ takes a moment frame and inserts extension joints %in the middle of each edge
as well as pinned anchors at each junction.

\section{The Homological Relations of Moment Frames}
\begin{figure}[h]
    \centering
        \includegraphics[width=\textwidth]{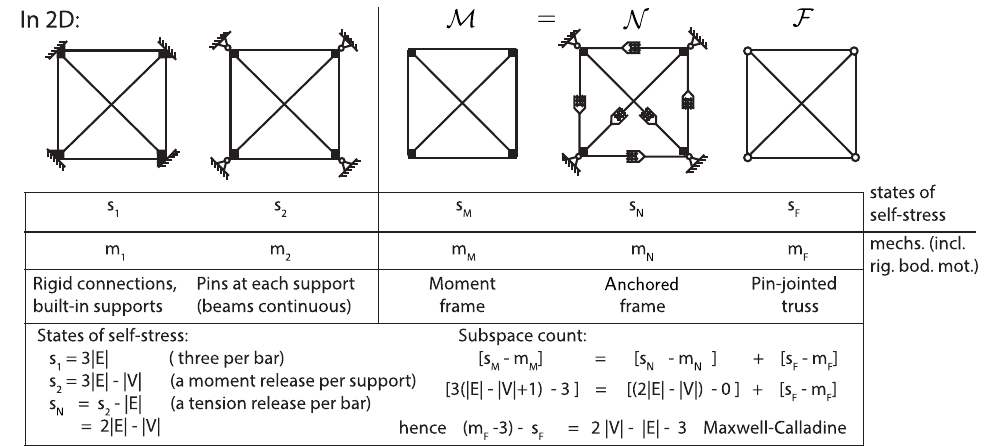}
            \caption{Five structures. The first has fully fixed joints and supports. The second releases the support moments. The final three are the fully-rigid frame, the anchored frame and the truss. The Maxwell-Calladine count for trusses is ($3\times$) the circuit count for frames minus the count for the anchored frame.} 
        \label{five_structures}
\end{figure}

We motivate the following technical analysis with a simple example.
\begin{example}[Frame counts]
    Consider the 2D frame whose bars are all moment-connected and whose joints are all fully encastr{\'e} supports. Releases of constraints lead variously to the anchored frame and the pin-jointed truss (see Fig.~\ref{five_structures}). The number of states of self-stress in the anchored frame may be readily determined by considering the first two frames  incrementally removing constraints. Following Equation~\eqref{eq:euler_moment}, the number of states of self-stress for the moment frame is 3 times the circuit rank, including three rigid body motions in 2D. Together, these counts redevelop the Maxwell-Calladine count for the 2D truss.
\end{example}

The homological algebra allows us to mathematically formalise these relationships. Every short exact sequence of cosheaves~(\ref{eq:force_moment_ses}) induces a long exact sequence of homology for all indices $i$:
\begin{equation}
\label{eq:les}
    \cdots \xrightarrow{\pi_*} H_{i+1} \caln \xrightarrow{\vartheta} H_i \calf \xrightarrow{\phi_*} H_i \calm \xrightarrow{\pi_*} H_i \caln \xrightarrow{\vartheta} H_{i-1} \calf \xrightarrow{\phi_*} \cdots
\end{equation}
where $\phi_*$ and $\pi_*$ are induced maps between homology spaces (which simply restrict and project onto homology spaces). The maps $\vartheta: H_{i+1} \caln \to H_i \calf$ are called {\it connecting homomorphisms} following from the {\it snake/zig-zag lemma} in algebraic topology~\cite{AlgebraicHatcher2002}. Exactness means that for every term in~\eqref{eq:les} the image of the incoming linear map is the kernel of the outgoing map. The derivation is too much to go into here but see~\cite{Part1Cooperband2023} for an exposition. The takeaway is that this is a powerful method from algebraic topology that describes linear relations between homology spaces (here self-stresses and DOF).

\begin{figure}[ht]
    \centering
    \begin{subfigure}[b]{0.98\textwidth}\centering
        \includegraphics[scale = 1.0]{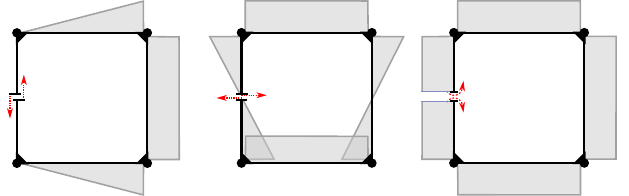}
        \caption{Three basis cycles for $H_1 \calm$ are pictured. Each is generated from a unit axial, shear or moment over the leftmost member. These pass to linearly independent elements of $H_1 \caln$.}
    \end{subfigure}
    
    \hfill
    \begin{subfigure}[b]{.58\textwidth}\centering
        \includegraphics[scale = 1.1]{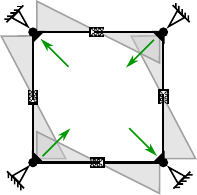}
        \caption{The final generator of $H_1 \caln$ with green resultant forces.}
    \end{subfigure}
    \hfill
    \begin{subfigure}[b]{0.4\textwidth}\centering
    \begin{tabular}{r |c c c}
         & $\calf$ & $\calm$ & $\caln$  \\
         \hline
         $s; \dim H_1$ & 0 & 3 & 4 \\
         $m; \dim H_0$ & 4 & 3 & 0
    \end{tabular}
    \caption{Table of cosheaf homology dimensions.}
    \end{subfigure}
    \hfill
    \vspace{1pt}
    \caption{Anchored cosheaf self-stresses over a box frame.}
    \label{fig:moment_square}
\end{figure}

Over a framework $(G, p)$ in $\R^2$ the homology space $H_0 \calf$ consists of pin-jointed kinematic degrees of freedom containing at least $3$ degrees of rigid body DOF. These trivial DOF are precisely those that generate $H_0 \calm$. The long exact sequence of the trio $\calf, \calm, \caln$ is

% As a result, the surjective homology map $\phi: H_0 \calf \to H_0 \calm$ maps rigid body motion to rigid body motion \textcolor{red}{That's true, but sounds a bit odd - trivial rigid body motion is tautological}.Preemptively removing these trivial degrees of freedom, 

\begin{equation}
\label{eq:moment_les}
    0 \to H_1 \calf \xrightarrow{\phi_*} H_1 \calm \xrightarrow{\pi_*} H_1 \caln \xrightarrow{\vartheta} \scrm \calf \oplus \scrr \calf \xrightarrow{\phi_*} \scrr \calm \to 0
\end{equation}
where the homology space $H_0 \caln$ must be zero because $\phi_*: H_0 \calf \to H_0 \calm$ is surjective. This homology map $\phi_*$ maps rigid body truss DOFs to rigid body frame DOFs with kernel $\scrm \calf$ (mechanisms aren't frame DOFs and so $\phi_* (\scrm \calf) \subset \im \bdd_\calm$). We can simplify this sequence by removing both rigid body DOF terms to form a reduced long exact sequence
\begin{equation}
\label{eq:moment_les_reduced}
    0 \to H_1 \calf \xrightarrow{\phi_*} H_1 \calm \xrightarrow{\pi_*} H_1 \caln \xrightarrow{\vartheta} \scrm \calf \to 0
\end{equation}
where the connecting homomorphism $\vartheta: H_1 \caln \to \scrm \calf$ is surjective onto mechanisms. Here $\vartheta$ is simply the linear map that inputs an anchor frame self-stress and outputs the resultant shear forces at nodes --- the green force arrows in Fig.~\ref{fig:moment_square} (b) and Fig.~\ref{fig:moment_desargues} (a). It is interesting to interpret the physical meaning of sequence~\eqref{eq:moment_les_reduced}:
\begin{itemize}
    \item At $H_1 \calf$, every axial self-stress is also a frame-self stress in $H_1 \calm$ ($\phi_*$ is injective).
    \item Every anchored cosheaf self-stress in $H_1 \caln$ is a combination of a moment frame self-stress in $H_1 \calm$, and a stress corresponding to a pin-jointed truss mechanism. Formally this is the statement $H_1 \caln \iso \im \pi_* \oplus \im \vartheta^\top$ where $\vartheta^\top$ is the adjoint/transpose of $\vartheta$.
    \item The last map $\vartheta: H_1 \caln \to \scrm \calf$ is surjective, so every pin-joint mechanism in $\scrm \calf$ is the image of an anchored frame self-stress $H_1 \caln$ with force resultants encoding the mechanism velocities.
\end{itemize}

\begin{figure}[ht]
    \centering
    \includegraphics[scale = 0.9]{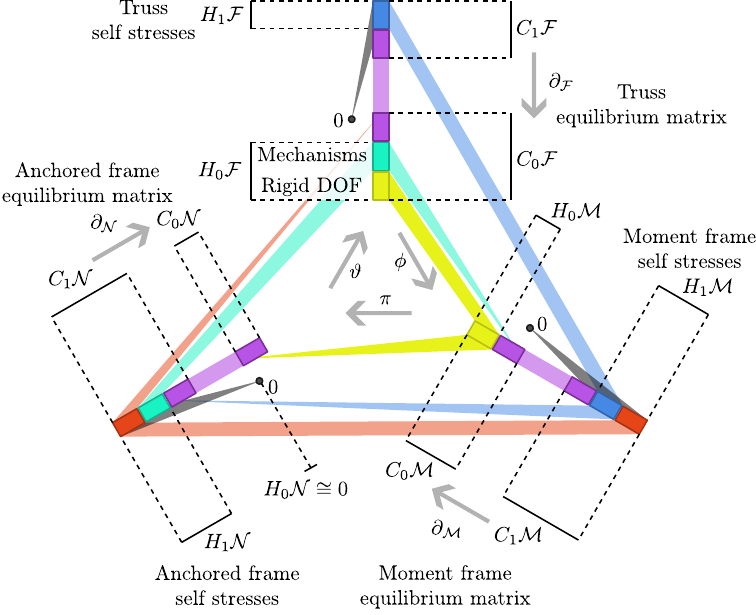}
    \caption{A schematic of the equilibrium operations in the three structural systems (truss, frame, and anchor frame) as well as the linear relations in the long exact sequence~\eqref{eq:moment_les} in clockwise order. The blue regions are pure axial self-stresses, the red regions are mixed frame self-stresses, the cyan regions are mechanisms/pure anchored frame self-stresses, and the yellow regions are rigid body DOF. The purple regions correspond to the dimensions of full rank of each equilibrium matrix.}
    \label{fig:linalg_les}
\end{figure}

These algebraic relations are sketched in Fig.~\ref{fig:linalg_les}. The last point is significant enough for its own statement.

\begin{theorem}[Mechanisms from the anchor-frame]
    Every pin-joint mechanism of a framework $(G, p)$ follows from the shear resultants of an anchored frame self-stress over $(G, p)$ (i.e. support reactions).
\end{theorem}

\begin{example}[Moments over a box frame]
\label{ex:moment_box_frame}
    Suppose that $(G, p)$ is the simple square frame embedded in $\R^2$. There are clearly no truss axial self-stresses so line~\eqref{eq:moment_les_reduced} reduces to a short exact sequence. By the isomorphism $H_1 \caln \iso H_1 \calm \oplus \scrm \calf$,  we see that there are $3+1 = 4$ states of anchor frame self-stress. The first three come from self-moments in $H_1 \calm$, generated by a unit axial force, shear force, and moment at a particular edge. These are pictured in Fig.~\ref{fig:moment_square} (a). The last state imparts forces which are absorbed by the anchors, pictured in Fig.~\ref{fig:moment_square} (b). These forces also inform a mechanism of a square pin-jointed truss, warping the square into a parallelogram.
\end{example}

\subsection{An Anchored Stress Counting Rule}

The Euler count~\eqref{eq:euler} of the exact sequence~\eqref{eq:moment_les_reduced} (as a chain complex) is the following alternating sum:
\begin{equation}
\label{eq:euler_moment_les}
    \mathfrak{X}\eqref{eq:moment_les_reduced} = \left( | H_1 \calf|- |\scrm \calf| \right) + |H_1 \caln| - |H_1 \calm| = 0
\end{equation}
which by exactness equals zero. These homology dimensions can be reformed in terms of cell numbers. The left-most term is the reduced Maxwell count~\eqref{eq:maxwell_rule} without rigid DOF while the right-most term $|H_1 \calm|$ has the cycle dimension $3(|E| - |V| +1)$ following from Equation~(\ref{eq:euler_moment}). Since $H_0 \caln=0$, then a simple Euler count gives $|H_1 \caln| = |C_1 \caln| - |C_0 \caln| = 2|E| - |V|$. Inserting these into Equation~\eqref{eq:euler_moment_les}, we have the following: 

\begin{theorem}[Anchor-frame stress count]
\label{thm:anchor_euler}
    Over a framework $(G, p)$, the dimension of anchored frame self-stresses $|H_1 \caln|$ is equal to a multiple of the cycle count minus the reduced Maxwell count~\eqref{eq:maxwell_rule}. Thus 
    %In $\R^2$ this is the equation:
    \begin{eqnarray}
    \mathrm{in} \  \R^2 : \quad    |H_1 \caln| & = & 3 (|E| - |V| +1) - (|E| - 2|V| + 3) = 2|E| -|V| \\
    %\end{equation}
    %In $\R^3$ this is the equation:
    %\begin{equation}
    \mathrm{in} \ \R^3 : \quad 
        |H_1 \caln|  & = &  6 (|E| - |V| +1) - (|E| - 3|V| + 6) = 5|E| - 3|V|
    \end{eqnarray}
\end{theorem}

To the authors' knowledge, this is the first time the Maxwell count has been combined with the graph-theoretic cycle count in an application. Visually in Fig.~\ref{fig:linalg_les}, Theorem~\ref{thm:anchor_euler} is the statement that the dimensions of the red and cyan regions can be counted two ways. The first is directly counting $|H_1 \caln|$. The second is through summing $|C_1 \calm| + |C_0 \calf|$ then subtracting the dimensions of $C_0 \calm$ and $C_1 \calf$, removing the blue, yellow and purple regions with only red and cyan regions remaining.

\begin{figure}[ht]
    \centering
    \begin{subfigure}[b]{0.68\textwidth}\centering
        \includegraphics[scale = 1.0]{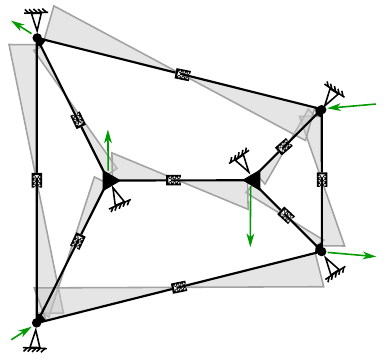}
        \caption{The generator of $H_1 \caln$ orthogonal to the frame self-stresses in $\pi_* (H_1 \calm)$. This self-stress of $\caln$ is mapped to a Desargues truss mechanism by $\vartheta$.}
    \end{subfigure}
    \hfill
    \begin{subfigure}[b]{0.3\textwidth}\centering
        \begin{tabular}{r |c c c}
             & $\calf$ & $\calm$ & $\caln$  \\
             \hline
             $s; \dim H_1$ & 1 & 12 & 12 \\
             $m; \dim H_0$ & 4 & 3 & 0
        \end{tabular}
    \caption{Table of cosheaf homology dimensions.}
    \end{subfigure}
    \vspace{1pt}
    \caption{An anchor-frame self-stress in $H_1 \caln$ is pictured over the Desargues configuration.}
    \label{fig:moment_desargues}
\end{figure}

\begin{example}[The Desargues anchor frame]
    Examine the Desargues configuration in Fig.~\ref{fig:moment_desargues} (a) with cosheaf homology dimensions listed in (b). The homology map $\phi_*: H_1 \calf \to H_1 \calm$ is injective and has rank 1, so by exactness $\ker\pi_* = \im \phi_*$ is dimension 1, and thus the rank-nullity theorem states that $\pi_*: H_1 \calm \to H_1 \caln$ has rank 11. The sole generator of $(\im \pi_*)^\perp \subset H_1 \caln$ is pictured in (a); other generators of $H_1 \caln$ are derived from ordinary frame self stresses $H_1 \calm$ of the Desargues frame.

    Moving down the long exact sequence~\eqref{eq:moment_les_reduced}, because $\vartheta: H_1 \caln \to \scrm \calf$ is surjective the anchor-frame stress in (a) generates the sole mechanism of the configuration. The residual shear forces when summed, drawn with green arrows, form the mechanism of the %pin-jointed Desargues
    truss (with anchors removed).

    From Theorem~\ref{thm:anchor_euler} we find that $|H_1 \caln| = 2|E| - |V| = 12$, in agreement with the table in Fig.~\ref{fig:moment_desargues} (b). The Maxwell count is $|H_1 \calf| - |H_0 \calf| = -3$ and the extended cycle count is $|H_1 \calm| - |H_0 \calm| = 9$. Taking the difference of these two counts is $9 - (-3) = 12$, the aforementioned count of anchored-frame self-stresses.

    The existence of the non-trivial self-stress and mechanism of the Desargues framework relies on the three vertical edges meeting at a point (in projective space), a delicate condition. %After perturbing the configuration this projective property will cease to hold.
    After perturbing the system both $H_1 \calf$ and $H_0 \calf$ decrease in dimension while the dimensions of $H_1 \calm$ and $H_1 \caln$ stay the same. What changes in sequence~\eqref{eq:moment_les_reduced} is the rank of the map $\phi_*: H_1 \calf \to H_1 \calm$ decreases, which by exactness means that $\pi_*: H_1 \calm \to H_1 \caln$ must increase in rank. Thus the cycle pictured in Fig.~\ref{fig:moment_desargues} (a) %, when perturbed,
    becomes a state of frame self-stress in $H_1 \calm$. This example suggests that the cosheaf $\caln$ can be used to better understand pin-jointed truss mechanics near singular points in its geometry.
\end{example}

% \newpage
% .
% \newpage
% \input{Zoe_main.tex}

\section{Conclusion}

We have described a new structural system, the anchored frame whose moment-dominated statics tie together the statics of pin-jointed trusses and rigid moment frames. The theory led to a counting rule for anchored frames incorporating both the Maxwell-Calladine rule and the graph cycle count. Moreover, each pin-jointed mechanism is encoded by the resultants of anchored frames. For future work, the authors intend to describe the discontinuous Airy stress functions~\cite{Williams2016} homologically.
% \textcolor{red}{Add another sentence to conclusion?}

\section*{Acknowledgments}
Zoe Cooperband was supported by the Air Force Office of Scientific Research. %under award number FA9550-21-1-033.
Bernd Schulze was partially supported by the ICMS Knowledge Exchange Catalyst Programme. The authors would like to thank Bill Baker, Chris Calladine, Robert Ghrist and Miguel Lopez for insightful discussions.

\printbibliography
\end{document}